\documentclass[reqno]{amsart}
\usepackage{amsmath, amsthm, amssymb, amstext}
\newtheorem{theorem}{Theorem}

\newtheorem{lemma}[theorem]{Lemma}

\newtheorem{corollary}[theorem]{Corollary}

\DeclareMathOperator*{\divergenz}{div}              %
\DeclareMathOperator*{\ints}{int}         %
\DeclareMathOperator*{\loc}{loc}         %
\DeclareMathOperator*{\essinf}{ess ~inf}         %
\DeclareMathOperator*{\esssup}{ess ~sup}         %

\def\ra{\rightarrow}

\def\cd{\cdot}

\def\N{\mathbb{N}}

\def\R{\mathbb{R}}
\def\C{\mathbb{C}}
\def\RN{\mathbb{R}^N}

\def\Om{\Omega}

\def\W{W^{1,p}(\Omega)}

\def\Lpx{L^{p(\cd)}(\Omega)}

\def\eps{\varepsilon}
\def\ph{\varphi}
\def\rand{\partial \Omega}
\def\into{\int_{\Omega}}

\def\W1p{W^{1,p}(\Omega)}
\def\Wx1p{W^{1,p(\cd)}(\Omega)}

\def\Linf{L^{\infty}(\Omega)}
\def\C1{\ints (C^1(\overline{\Om})_+)}

\numberwithin{theorem}{section}
\numberwithin{equation}{section}

\def\cprime{$'$}

\title[A priori bounds for elliptic equations with nonstandard growth]
      {A priori bounds for weak solutions to elliptic equations with nonstandard growth}
\author[P. Winkert]{Patrick Winkert}
\address{Technische Universit\"{a}t Berlin, Institut f\"{u}r Mathematik,\\ Stra\ss e des 17. Juni 136, 10623 Berlin, Germany}
\email{winkert@math.tu-berlin.de}
\author[R. Zacher]{Rico Zacher}
\address{Martin-Luther-Universit\"{a}t Halle-Wittenberg, Institut f\"{u}r Mathematik, Theodor-Lieser-Stra\ss e 5, 06120 Halle, Germany}
\email{rico.zacher@mathematik.uni-halle.de}

\subjclass[2010]{35J60, 35B45, 35J25}
\keywords{A priori estimates, De Giorgi iteration, Elliptic equations, Nonstandard growth, Partition of unity, Variable exponent spaces}
\begin{document}
\maketitle

\begin{abstract}
	In this paper we study elliptic equations with a nonlinear conormal derivative boundary condition involving nonstandard growth terms. By means of 
    the localization method and De Giorgi's iteration technique we derive global a priori bounds for weak solutions of such problems.
\end{abstract}

\section{Introduction}\label{section1}
The present paper is concerned with global a priori bounds for
elliptic equations with nonlinear conormal derivative boundary
conditions which may contain nonlinearities with variable growth
exponents. More precisely, let $\Om$ be a bounded domain in $\RN, N > 1$,
with Lipschitz boundary $\Gamma :=\rand$ and let $p \in
C(\overline{\Om})$ be a function that satisfies
$1<p^-:=\inf_{\overline{\Om}} p(x)$. We deal with elliptic equations
of the form
\begin{align}\label{problem}
    \begin{aligned}
        -\divergenz \mathcal{A} (x,u,\nabla u) & = \mathcal{B}(x,u,\nabla u)  & \hspace*{0.5cm} & \text{ in } \Om, \\
    \mathcal{A} (x,u,\nabla u)\cdot \nu & = \mathcal{C}(x,u)  && \text{ on } \Gamma,
    \end{aligned}
\end{align}
where $\nu(x)$ denotes the outer unit normal of $\Omega$ at $x\in
\Gamma$, and $\mathcal{A}, \mathcal{B}$ and $\mathcal{C}$ satisfy
suitable $p(x)$-structure conditions, see (H) below.

An important special case of (\ref{problem}) which fits in our
setting is given by
\begin{align*}
    -\Delta_{p(x)}u= \mathcal{B}(x,u,\nabla u) \quad \text{ in }  \Om,\qquad
      |\nabla
u|^{p(x)-2}\partial_\nu u= \mathcal{C}(x,u) \quad  \text{ on }
    \Gamma.
\end{align*}
Here the operator $\divergenz \mathcal{A}$ becomes the so-called
$p(x)$-Laplacian
\[
\Delta_{p(x)}u=\divergenz(|\nabla u|^{p(x)-2}\nabla u), \] which
reduces to the standard $p$-Laplacian if $p(x)\equiv p$.

In recent years there has been a growing interest in the study of
elliptic problems with a $p(x)$-structure, which are also termed
problems with nonstandard growth conditions. Equations of this type
appear in the study of non-Newtonian fluids with thermo-convective
effects (see Antontsev and Rodrigues \cite{2006-Antontsev}, Zhikov
\cite{1997-Zhikov}), electro-rheological fluids (see Diening
\cite{2002-Diening}, Rajagopal and R\r u\v zi\v cka
\cite{2001-Rajagopal}, R\r u\v zi\v cka \cite{2000-Ruzicka}), the
thermistor problem (see Zhikov \cite{1997-Zhikov2}), or the problem
of image recovery (see Chen et al. \cite{2006-Chen}).

Throughout the paper we impose the following conditions.
\begin{enumerate}
    \item[(H)]
    The functions $\mathcal{A}: \Om \times \R \times \RN \ra \RN$, $\mathcal{B}: \Om \times \R \times \RN \ra \R$, and
    $\mathcal{C}: \Gamma \times \R \ra \R$ are Carath\'eodory
    functions satisfying the subsequent structure
    conditions:
    \begin{align*}
        \text{(H1) \quad } & |\mathcal{A}(x,s,\xi)| \leq a_0|\xi|^{p(x)-1}+a_1|s|^{q_0(x)\frac{p(x)-1}{p(x)}}+a_2, && \text{ for a.a. } x\in \Om,\\
        \text{(H2) \quad } & \mathcal{A}(x,s,\xi) \cd \xi \geq a_3|\xi|^{p(x)}-a_4|s|^{q_0(x)}-a_5, && \text{ for a.a. } x\in \Om,\\
        \text{(H3) \quad } & |\mathcal{B}(x,s,\xi)| \leq b_0|\xi|^{p(x)\frac{q_0(x)-1}{q_0(x)}}+b_1|s|^{q_0(x)-1}+b_2, &&\text{ for a.a. } x\in \Om,\\
        \text{(H4) \quad } & |\mathcal{C}(x,s)| \leq c_0|s|^{q_1(x)-1}+c_1, && \text{ for a.a. } x\in \Gamma,
    \end{align*}
and for all $s \in \R$, and all $\xi \in \RN$. Here $a_i, b_j$ and
$c_l$ are positive constants, $p\in C(\overline{\Om})$ with
$\inf_{\overline{\Om}} p(x)>1$, and $q_0 \in C(\overline{\Om})$ as
well as $q_1\in C(\Gamma)$ are chosen such that
\begin{align*}
p(x) \leq q_0(x)<p^*(x),  \;\; x\in \overline{\Om}, \quad
\text{and}\quad p(x)\leq q_1(x)<p_*(x),\;\;x\in \Gamma,
\end{align*}
with the
critical exponents
\begin{align*}
        p^*(x)=
        \begin{cases}
            \frac{Np(x)}{N-p(x)} \quad & \text{ if } p(x) <N, \\
            +\infty & \text{ if } p(x) \geq N,
        \end{cases} \qquad
        p_*(x)=
        \begin{cases}
            \frac{(N-1)p(x)}{N-p(x)} \quad & \text{ if } p(x) <N, \\
            +\infty & \text{ if } p(x) \geq N.
        \end{cases}
\end{align*}
\end{enumerate}

A function $u \in \Wx1p$ is said to be a {\bf weak solution}
({\bf subsolution, supersolution}) of equation (\ref{problem}) if
\begin{align}\label{solution}
        \into \mathcal{A}(x,u,\nabla u) \cdot \nabla \ph dx =\,(\le,\,\ge)
        \into \mathcal{B}(x,u,\nabla u) \ph dx + \int_{\Gamma} \mathcal{C}(x,u) \ph d \sigma,
\end{align}
holds for all nonnegative test functions $\ph \in \Wx1p$, where $d
\sigma$ denotes the usual $(N-1)$-dimensional surface measure.

This definition makes sense, since thanks to assumption (H) the
integrals in (\ref{solution}) are finite, by H\"older's inequality
and embedding results for $\Wx1p$-functions, see below.

The main goal of this paper is to prove a priori bounds for weak
sub- and supersolutions, in particular for weak solutions of problem
(\ref{problem}). Using the notation $y_+=\max(y,0)$, our main result
reads as follows.
\begin{theorem}\label{maintheorem}
    Let the assumptions in (H) be satisfied. Then there exist
    positive constants $\alpha=\alpha(p,q_0,q_1)$ and
    $C=C(p,q_0,q_1,a_3,a_4,a_5,b_0,b_1,b_2,c_0,c_1,N,\Omega)$ such
    that the following assertions hold.
    \begin{enumerate}
        \item[(i)]
            If $u \in \Wx1p$ is a weak subsolution of (\ref{problem}) then
            \begin{align*}
                \esssup_{\Om} u \leq 2 \max\left(1, C
                \left [\into u_+^{q_0(x)} dx + \int_{\Gamma} u_+^{q_1(x)} d \sigma \right
                ]^\alpha\right).
            \end{align*}
        \item[(ii)]
            If $u \in \Wx1p$ is a weak supersolution of (\ref{problem}) then
            \begin{align*}
                \essinf_{\Om} u \geq -2 \max\left(1, C
                \left [\into (-u)_+^{q_0(x)} dx + \int_{\Gamma} (-u)_+^{q_1(x)} d \sigma \right ]^\alpha\right).
            \end{align*}
    \end{enumerate}
\end{theorem}
Note that the constants $a_0,a_1,a_2$, which appear in (H1), do not
play any role in determining the constants $\alpha$ and $C$. The
finiteness of the right-hand sides in (i) and (ii) is a consequence
of the compact embedding $\Wx1p \ra L^{q_0(\cd)}(\Om)$ and the fact
that the trace operator is a bounded operator from $\Wx1p$ into
$L^{q_1(\cd)}(\Gamma)$ (see Fan et al. \cite[Theorem 1.3]{2001-Fan}
and Fan \cite[Corollary 2.4]{2008-Fan}). We further point out that
we merely assume continuity for the variable exponents $p,q_0,$ and
$q_1$; log-H\"older continuity conditions are not required.

Our proof of Theorem \ref{maintheorem} uses De Giorgi's iteration
technique and the localization method. By means of the latter we are
able to reduce the estimates involving variable exponents to ones
with constant exponents, which then also allows us to apply
classical embedding results. This crucial step in the proof is
achieved by means of an appropriate partition of unity.

By the definition of sub- and supersolution of (\ref{problem}) one
easily verifies that a weak solution is both a weak subsolution and
a weak supersolution. Hence we have the following.
\begin{corollary}
    Let the assumptions (H) be satisfied and let $u \in \Wx1p$ be a weak solution of (\ref{problem}).
    Then $u\in L^\infty(\Omega)$ and the estimates in (i) and (ii)
    from Theorem \ref{maintheorem} are valid.
\end{corollary}
The main novelty of the paper consists in the generality of the
assumptions needed to establish the boundedness of weak solutions to
(\ref{problem}). In particular the assumptions on the nonlinearity
$\mathcal{C}$ are rather general, allowing for a growth term with
variable exponent, which seems to be optimal. Another novelty is the
use of the localization technique in the context of global a priori
estimates for problems with variable exponents and nonlinear
conormal derivative boundary conditions.

Let us comment on some relevant known results on elliptic problems
with $p(x)$-structure. Local boundedness of solutions to the
equation
\begin{align}\label{fan}
    -\divergenz \mathcal{A} (x,u,\nabla u) & = \mathcal{B}(x,u,\nabla u)  \quad \text{ in } \Om,
\end{align}
has been studied by Fan and Zhao \cite{1999-Fan}. There it is shown that
under suitable structure conditions every weak solution $u$ of
(\ref{fan}) (corresponding to test functions $\ph \in
W^{1,p(x)}_0(\Om)$) belongs to $L^{\infty}_{\loc}(\Om)$, and if in
addition $u$ is bounded on the boundary $\Gamma$, then $u \in
\Linf$. The proof uses De Giorgi iterations as well. Recently, Gasi{\'n}ski and Papageorgiou (see \cite[Proposition 3.1]{2011-Gasinski}) studied global a priori bounds for weak solutions to the equation
\begin{align}\label{gasinski}
    \begin{aligned}
        -\Delta_{p(x)} u & = g(x,u)  & \hspace*{0.5cm} & \text{ in } \Om, \\
    \frac{\partial u}{\partial \nu} & = 0  && \text{ on } \Gamma,
    \end{aligned}
\end{align}
where the Carath\'{e}odory function $g: \Om \times \R \ra \R$
satisfies a subcritical growth condition. They proved that every
weak solution $u \in \Wx1p$ of problem (\ref{gasinski}) belongs to
$\Linf$ provided $p \in C^1(\overline{\Omega})$ satisfying
$1<\min_{x \in \overline{\Omega}}p(x)$.

$L^{\infty}$-estimates for solutions of (\ref{problem}) in case $p(x)\equiv p$ with $q_0(x)=q_1(x)\equiv p$ have been established by the first author in \cite{2010-Winkert1,2010-Winkert2} following Moser's iteration technique (for constant $p$ see also Pucci and Servadei \cite{2008-Pucci}).

Concerning boundedness and regularity results for problems of type (\ref{fan}), in particular for the special case
\begin{align*}
    -\divergenz (|\nabla u|^{p(x)-2}\nabla u) =0,
\end{align*}
we further refer to Acerbi and Mingione
\cite{2001-Acerbi,2002-Acerbi}, Antontsev and Consiglieri
\cite{2009-Antontsev}, Chiad\`{o} Piat and Coscia
\cite{1997-ChiadoPiat}, Diening et al. \cite{2007-Diening}, Eleuteri
and Habermann \cite{2008-Eleuteri}, Fan \cite{2007-Fan,2010-Fan},
Fan and Zhao \cite{2000-Fan}, Habermann and Zatorska-Goldstein
\cite{2008-Habermann}, Harjulehto et al. \cite{2007-Harjulehto},
Liskevich and Skrypnik \cite{2010-Liskevich}, Lukkari
\cite{2010-Lukkari,2009-Lukkari} and the references given therein.

The paper is organized as follows. In Section 2 we fix some notation
and recall the definition of the variable exponent spaces $\Lpx$ and
$\Wx1p$. We further state a lemma on sequences of numbers which will
be needed for the De Giorgi iterations. The main result is proved in
Section 3. From the structure conditions we first derive truncated
energy estimates. These are then used, together with the
localization method and embedding results, to prove suitable
iterative inequalities, which in turn imply the desired a priori
bounds.
\section{Notations and preliminaries}\label{section2}
Suppose that $\Om$ is a bounded domain in $\RN$ with Lipschitz
boundary $\Gamma$ and let $p \in C(\overline{\Om})$ with $p(x)>1$
for all $x \in \overline{\Om}$. We set $p^-:=\min_{x \in
\overline{\Om}} p(x)$ and $p^+:=\max_{x \in \overline{\Om}} p(x)$,
then $p^->1$ and $p^+<\infty$. By $\Lpx$ we identify the variable
exponent Lebesgue space which is defined by
\begin{align*}
    \Lpx = \left \{u ~ \Big | ~ u: \Om \ra \R \text{ is measurable and } \into |u|^{p(x)}dx< +\infty \right \}
\end{align*}
equipped with the Luxemburg norm
\begin{align*}
    \|u\|_{\Lpx} = \inf \left \{ \tau >0 : \into \left |\frac{u(x)}{\tau} \right |^{p(x)}dx \leq 1  \right \}.
\end{align*}
The variable exponent Sobolev space $\Wx1p$ is defined by
\begin{align*}
    \Wx1p= \{u \in \Lpx : |\nabla u| \in \Lpx \}
\end{align*}
with the norm
\begin{align*}
    \|u\|_{\Wx1p}= \|\nabla u \|_{\Lpx}+\|u\|_{\Lpx}.
\end{align*}
For more information and basic properties of variable exponent
spaces we refer the reader to the papers of Fan and Zhao
\cite{2001-Fan2}, Kov{\'a}{\v{c}}ik and R{\'a}kosn{\'{\i}}k
\cite{1991-Kovacik} and the recent monograph of Diening et al. \cite{2011-Diening}. If $p(x)\equiv p$ is a constant, the usual
Sobolev space $\W1p$ is endowed with the norm
\begin{align*}
    \|u\|_{\W1p} = \left ( \into |\nabla u|^p dx +\into |u|^p dx \right )^{\frac{1}{p}}.
\end{align*}
For $q_0 \in C(\overline{\Om})$ and $q_1\in C(\Gamma)$ (as in (H))
we define
 \begin{align*}
        \begin{split}
            & q_0^+=\max_{\overline{\Om}} q_0(x), \quad q_0^-=\min_{\overline{\Om}}q_0(x),\\
            & q_1^+=\max_{\Gamma} q_1(x), \quad q_1^-=\min_{\Gamma}q_1(x).
        \end{split}
\end{align*}
For $s\in [1,\infty)$ we further use the notation
\begin{align*}
        s^*=
        \begin{cases}
            \frac{Ns}{N-s} \quad & \text{ if } s <N, \\
            +\infty & \text{ if } s \geq N,
        \end{cases} \qquad
        s_*=
        \begin{cases}
            \frac{(N-1)s}{N-s} \quad & \text{ if } s <N, \\
            +\infty & \text{ if } s \geq N.
        \end{cases}
\end{align*}

%

The following lemma concerning the geometric convergence of
sequences of numbers will be needed for the De Giorgi iteration
arguments below. It can be found, for example in
\cite{2010-Vergara}. The case $\delta_1=\delta_2$ is contained in
\cite[Chapter II, Lemma 5.6]{1968-Lady}, see also \cite[Chapter I,
Lemma 4.1]{1993-DiBenedetto}.
\begin{lemma}\label{lemma4}
    Let $\{Y_n\}, n=0,1,2,\ldots,$ be a sequence of positive numbers, satisfying the recursion inequality
    \begin{align*}
        Y_{n+1} \leq K b^n \left (Y_n^{1+\delta_1}+ Y_n^{1+\delta_2} \right ) , \quad n=0,1,2, \ldots,
    \end{align*}
    for some $b>1,\,K>0$, and $\delta_2\geq \delta_1>0$. If
    \begin{align*}
        Y_0 \leq (2K)^{-\frac{1}{\delta_1}} b^{-\frac{1}{\delta_1^2}},
    \end{align*}
    then
    \begin{align*}
        Y_n \leq (2K)^{-\frac{1}{\delta_1}} b^{-\frac{1}{\delta_1^2}} b^{-\frac{n}{\delta_1}}, \quad n \in \N,
    \end{align*}
    in particular $\{Y_n\} \ra 0$ as $n \ra \infty$.
\end{lemma}
\section{Truncated energy estimates and proof of Theorem \ref{maintheorem}}
Our proof of the sup-bounds for weak subsolutions of (\ref{problem})
is based on the following lemma on truncated energy estimates.
\begin{lemma}\label{lemma1}
        Let the conditions in (H) be satisfied. If $u$ is a weak subsolution of (\ref{problem}), there holds
        \begin{align*}
            \int_{A_k} |\nabla u|^{p(x)} dx  \leq d_1 \int_{A_k} u^{q_0(x)}dx+d_2 \int_{\Gamma_k} u^{q_1(x)}d \sigma,
        \end{align*}
    where
    \begin{align*}
        A_k=\{x \in \Om: u(x)>k\}, \qquad \Gamma_k = \{x \in \Gamma: u(x)>k\},  \quad k \geq 1,
    \end{align*}
    and $d_1=2a_3^{-1}(a_4+a_5+b_1+b_2+b_0 \eps^{-(q_0^++1)})$, $d_2=2a_3^{-1}(c_0+c_1)$, and $\eps=\min(1,\frac{a_3}{2b_0})$.
\end{lemma}

\begin{proof}
    Let $u \in \Wx1p$ be a weak subsolution of (\ref{problem}) and let $k \geq 1$. Taking $\ph=(u-k)_+=\max(u-k,0) \in \Wx1p$
    (see \cite[Lemma 3.2]{2009-Le}) as test function in
    (\ref{solution}) with the '$\le\,$'-sign
    we obtain
    \begin{align}\label{3.1}
        \begin{split}
            & \int_{A_{k}} \mathcal{A}(x,u,\nabla u) \cdot \nabla (u-k) dx\\
            & \leq \int_{A_{k}} \mathcal{B}(x,u,\nabla u) (u-k)dx + \int_{\Gamma_k} \mathcal{C}(x,u) (u-k) d \sigma.
        \end{split}
    \end{align}
    Using the structure condition (H2) we estimate the left-hand side of
    (\ref{3.1}) as follows.
    \begin{align}\label{3.2}
        \begin{split}
            & \int_{A_k} \mathcal{A}(x,u,\nabla u ) \cdot \nabla (u-k) dx \\
            & = \int_{A_k} \mathcal{A}(x,u,\nabla u ) \cdot \nabla u dx \\
            & \geq \int_{A_k} \left (a_3 |\nabla u|^{p(x)}-a_4|u|^{q_0(x)}-a_5 \right ) dx\\
            & \geq a_3 \int_{A_k} |\nabla u|^{p(x)} dx - (a_4+a_5) \int_{A_k} |u|^{q_0(x)} dx,
         \end{split}
    \end{align}
    as $u^{q_0(x)}>u>1$ in $A_k$.
    Now, we are going to estimate the right-hand side of (\ref{3.1}). By Young's inequality with $\eps \in
    (0,1]$ and condition (H3) we have
    \begin{align}\label{3.3}
        \begin{split}
            & \int_{A_{k}} \mathcal{B}(x,u,\nabla u) (u-k)dx \\
            & \leq \int_{A_{k}} \left [b_0 |\nabla u|^{p(x)\frac{q_0(x)-1}{q_0(x)}}+b_1 |u|^{q_0(x)-1}+b_2 \right ]  (u-k) dx\\
            & \leq b_0 \int_{A_{k}} \left [ \eps^{\frac{q_0(x)-1}{q_0(x)}} |\nabla u|^{p(x)\frac{q_0(x)-1}{q_0(x)}} \eps^{-\frac{q_0(x)-1}{q_0(x)}}u \right ] dx+(b_1+b_2) \int_{A_k}|u|^{q_0(x)}dx\\
            & \leq b_0 \int_{A_k} \eps |\nabla u|^{p(x)} dx + b_0 \int_{A_k} \eps^{-(q_0(x)-1)} u^{q_0(x)}dx +(b_1+b_2) \int_{A_k}|u|^{q_0(x)}dx \\
            & \leq \eps b_0 \int_{A_k}|\nabla u|^{p(x)} dx + \left ( b_0 \eps^{-(q_0^+-1)}+b_1+b_2\right ) \int_{A_k} u^{q_0(x)} dx.
        \end{split}
    \end{align}
    Thanks to condition (H4), the boundary integral can be estimated
    through
    \begin{align}\label{3.4}
        \begin{split}
            \int_{\Gamma_k} \mathcal{C}(x,u) (u-k) d \sigma
                & \leq \int_{\Gamma_k} (c_0 |u|^{q_1(x)-1}+c_1)(u-k) d \sigma \\
            & \leq (c_0+c_1) \int_{\Gamma_k} u^{q_1(x)} d \sigma,
        \end{split}
    \end{align}
    as $u>1$ on $\Gamma_k$. Combining (\ref{3.1})--(\ref{3.4}) and choosing $\eps = \min(1,\frac{a_3}{2b_0})$
    gives
    \begin{align*}
        & \frac{a_3}{2}\int_{A_k} |\nabla u|^{p(x)} dx\\
        &\leq \left (a_4+a_5+b_1+b_2+b_0 \eps^{-(q_0^++1)} \right) \int_{A_k} u^{q_0(x)}dx+(c_0+c_1) \int_{\Gamma_k} u^{q_1(x)}d \sigma.
    \end{align*}
    Dividing the last inequality by $\frac{a_3}{2} >0$ yields the assertion of the lemma.
\end{proof}

The corresponding result for supersolutions reads as follows.

\begin{lemma}\label{lemma2}
        Let the conditions in (H) be satisfied. If $u$ is a weak supersolution of (\ref{problem}), there holds
        \begin{align*}
            \int_{\tilde{A}_k} |\nabla u|^{p(x)} dx  \leq d_1 \int_{\tilde{A}_k} (-u)^{q_0(x)}dx+d_2
            \int_{\tilde{\Gamma}_k} (-u)^{q_1(x)}d \sigma ,
        \end{align*}
    where
    \begin{align*}
        \tilde{A}_k=\{x \in \Om: -u(x)>k\}, \qquad \tilde{\Gamma}_k = \{x \in \Gamma: -u(x)>k\},    \quad k \geq 1,
    \end{align*}
    and $d_1$ and $d_2$ are the same constants as in Lemma \ref{lemma1}.
\end{lemma}

\begin{proof}
    The proof is analogous to the previous one. We take $\ph=-(u+k)_-=-\min(u+k,0) \geq 0$ as test
    function in (\ref{solution}), which now holds with the
    '$\ge\,$'-sign,
    and use the same arguments as in the proof of Lemma \ref{lemma1}. This yields the asserted
    inequality.
\end{proof}

Now we are in position to prove the main result of this paper.

\begin{proof}[Proof of Theorem \ref{maintheorem}]

${}$

{\bf (i) Definition of the iteration variables $Z_n$, $\tilde{Z}_n$,
and basic estimates.} Let now
    \begin{align*}
        k_n=k \left (2 -  \frac{1}{2^{n}} \right ), \quad n=0,1,2, \ldots ,
    \end{align*}
    with $k \geq 1$ specified later and put
    \begin{align*}
        Z_n:=\int_{A_{k_n}} (u-k_n)^{q_0(x)}dx, \qquad \tilde{Z}_n:= \int_{\Gamma_{k_n}} (u-k_n)^{q_1(x)} d \sigma.
    \end{align*}
        We have
    \begin{align*}
        \begin{split}
            Z_n
            & \geq \int_{A_{k_{n+1}}} (u-k_n)^{q_0(x)} dx \geq \int_{A_{k_{n+1}}} u^{q_0(x)} \left (1-\frac{k_n}{k_{n+1}} \right )^{q_0(x)}dx\\
            & \geq \int_{A_{k_{n+1}}} \frac{1}{2^{q_0(x) (n+2)}} u^{q_0(x)}dx,
        \end{split}
    \end{align*}
    and thus
    \begin{align}\label{3.5}
        \int_{A_{k_{n+1}}} u^{q_0(x)} dx \leq 2^{q_0^+(n+2)} Z_n.
    \end{align}
    Analogously, we see that
    \begin{align}\label{3.6}
        \int_{\Gamma_{k_{n+1}}} u^{q_1(x)} d \sigma \leq 2^{q_1^+(n+2)}
        \tilde{Z}_n.
    \end{align}
From (\ref{3.5})--(\ref{3.6}) and Lemma \ref{lemma1} with $k$ being
replaced by $k_{n+1} \geq 1$ it follows that
    \begin{align}\label{3.8}
            \int_{A_{k_{n+1}}} |\nabla (u-k_{n+1})|^{p(x)} dx  \leq d_3 a^n (Z_n + \tilde{Z}_n),
    \end{align}
    where $d_3=\max\left (d_1 2^{2q_0^+},d_2 2^{2q_1^+}\right )$ and $a=\max \left (2^{q_0^+}, 2^{q_1^+} \right)$.

Furthermore, we have
    \begin{align}\label{3.7}
        \begin{split}
            |A_{k_{n+1}}|
            & \leq \int_{A_{k_{n+1}}} \left (\frac{u-k_n}{k_{n+1}-k_n}\right )^{q_0(x)} dx\\
            & \leq \int_{A_{k_{n}}} \frac{2^{q_0(x)(n+1)}}{k^{q_0(x)}} (u-k_n)^{q_0(x)}dx \\
            & \leq \frac{2^{q_0^+(n+1)}}{k^{q_0^-}} \int_{A_{k_{n}}} (u-k_n)^{q_0(x)}dx \\
            & = \frac{2^{q_0^+(n+1)}}{k^{q_0^-}} Z_n.
        \end{split}
    \end{align}

${}$

{\bf (ii) Partition of unity.} By compactness of $\overline{\Om}$,
for any $R>0$ there exists a finite open cover
$\{B_i(R)\}_{i=1,\ldots, m}$ of balls $B_i:=B_i(R)$ with radius $R$
such that $\overline{\Om}\subset \bigcup_{i=1}^m B_i(R)$. Moreover,
since $p\in C(\overline{\Om})$, $q_0 \in C(\overline{\Om})$, and
$q_1\in C(\Gamma)$, these functions are uniformly continuous on
$\overline{\Om}$ and $\Gamma$, respectively. Recalling that
\begin{align*}
p(x) \leq q_0(x)<p^*(x),  \;\; x\in \overline{\Om}, \quad
\text{and}\quad p(x)\leq q_1(x)<p_*(x),\;\;x\in \Gamma,
\end{align*}
we may take $R>0$ small enough such that
\begin{align*}
        p^+_i \leq q_{0,i}^+ <(p^-_i)^*, \qquad p^+_i \leq q_{1,i}^+<(p^-_i)_{*}, \qquad i=1, \ldots , m,
    \end{align*}
 where
    \begin{align*}
        \begin{split}
            & p^+_i=\max_{B_i\cap \overline{\Om}} p(x), \quad q_{0,i}^+=\max_{B_i\cap \overline{\Om}}q_0(x),\\
            & p^-_i=\min_{B_i\cap \overline{\Om}} p(x), \quad q_{1,i}^+=\max_{B_i\cap \Gamma}
            q_1(x).
        \end{split}
    \end{align*}
We next choose a partition of unity $\{\xi_i\}_{i=1}^m \subset
C^{\infty}_0 (\RN)$ associated to the open cover
$\{B_i(R)\}_{i=1,\ldots, m}$ (see e.g.\ \cite[Thm.
6.20]{1973-Rudin}), that is, we have
    \begin{align*}
        \mbox{supp}\,\xi_i\subset B_i,\;\;\; 0 \leq \xi_i \leq 1,\;\;i=1,\ldots,m,\quad \mbox{and}\;\;\sum_{i=1}^m \xi_i=1\;\;\mbox{on}\;
        \overline{\Omega}.
    \end{align*}
Let $L>0$ be a positive constant such that
\[
|\nabla \xi_i| \leq L,\quad i=1,\ldots,m.
\]

${}$

{\bf (iii) Estimating the gradient term in (\ref{3.8}) from below.}
Using the partition of unity from step (ii) we have
    \begin{align}\label{3.9}
        \begin{split}
            & \int_{A_{k_{n+1}}} |\nabla (u-k_{n+1})|^{p(x)}dx \\
            & =  \int_{A_{k_{n+1}}} |\nabla (u-k_{n+1})|^{p(x)} \sum_{i=1}^m \xi_i dx\\
            & \geq  \sum_{i=1}^m  \int_{A_{k_{n+1}}} ( |\nabla (u-k_{n+1})|^{p^-_i}-1)  \xi_i dx\\
            & \geq \sum_{i=1}^m \int_{A_{k_{n+1}}} |\nabla (u-k_{n+1})|^{p^-_i} \xi_i^{p^-_i} dx - m |A_{k_{n+1}}|,
        \end{split}
    \end{align}
    as $\xi_i  \geq \xi_i^{p^-_i}$. From (\ref{3.9}) we trivially deduce that for all $i=1, \ldots,
    m$,
    \begin{align}\label{3.10}
        \begin{split}
            \int_{A_{k_{n+1}}} |\nabla (u-k_{n+1})|^{p(x)}dx \geq  \int_{A_{k_{n+1}}} |\nabla (u-k_{n+1})|^{p^-_i} \xi_i^{p^-_i} dx - m |A_{k_{n+1}}|.
        \end{split}
    \end{align}
    Combining (\ref{3.8}) and (\ref{3.10}) and using (\ref{3.7}) yields
    \begin{align}\label{3.11}
        \begin{split}
             & \int_{A_{k_{n+1}}} |\nabla (u-k_{n+1})|^{p^-_i} \xi_i^{p^-_i} dx \leq d_4 a^n (Z_n + \tilde{Z}_n)
        \end{split}
    \end{align}
    for any $i=1,\ldots,m$, with the positive constant $d_4=d_3+m2^{q_0^+}$.

${}$

{\bf (iv) Estimating $Z_{n+1}$.} Next we want to derive a suitable
estimate for the term $Z_{n+1}$ from above. To this end, we make
again use of the partition of unity introduced in step (ii). We have
    \begin{align}\label{3.12}
        \begin{split}
            Z_{n+1} &=\int_{A_{k_{n+1}}} (u-k_{n+1})^{q_0(x)}dx\\
            & =\int_{A_{k_{n+1}}} (u-k_{n+1})^{q_0(x)} \left ( \sum_{i=1}^m \xi_i\right )^{q_0^+}dx \\
            & \leq \int_{A_{k_{n+1}}} (u-k_{n+1})^{q_0(x)} m^{q_0^+} \sum_{i=1}^m \xi_i^{q_0^+} dx \\
            & \leq m^{q_0^+} \sum_{i=1}^m \int_{A_{k_{n+1}}} (u-k_{n+1})^{q_0(x)} \xi_i^{q_0(x)} dx \\
            & \leq m^{q_0^+} \sum_{i=1}^m \left [\int_{A_{k_{n+1}}}\!\!\!\! (u-k_{n+1})^{q_{0,i}^+}\xi_i^{q_{0,i}^+}dx +
            \int_{A_{k_{n+1}}}\!\!\!\! (u-k_{n+1})^{q_{0,i}^-}\xi_i^{q_{0,i}^-} dx \right ],
        \end{split}
    \end{align}
where we have set $q_{0,i}^-=\min_{B_i\cap \overline{\Om}}q_0(x)$.
Observe that $p^-_i \leq q^-_{0,i} \leq q^+_{0,i} <(p^-_i)^*$ for
all $i=1, \ldots, m$.

Let now $i\in \{1,\ldots,m\}$ be fixed, and suppose that $r\in
\{q^-_{0,i},q^+_{0,i}\}$. Then $p^-_i \leq r <(p^-_i)^*$ and $r \leq
q^+$, where $q^+=\max(q^+_0,q^+_1)$. By H\"{o}lder's inequality and
the continuous embedding $W^{1,p^-_i}(\Om) \hookrightarrow
L^{(p^-_i)^*}(\Om)$, we may estimate as follows.
    \begin{align}\label{3.13}
        \begin{split}
            & \into (u-k_{n+1})_+^{r}\xi_i^{r}dx\\
            & \leq \left (\into (u-k_{n+1})_+^{(p^-_i)^*} \xi_i^{(p^-_i)^*}dx \right)^{\frac{r}{(p^-_i)^*}}
            |A_{k_{n+1}}|^{1-\frac{r}{(p^-_i)^*}}  \\
            & \leq  C^{q^+} \left ( \into  \left [|\nabla  [(u-k_{n+1})_+ \xi_i]|^{p^-_i}\right. \right.\\
            & \qquad \left.\left. +|(u-k_{n+1})_+ \xi_i|^{p^-_i} \right ] dx \right)^{\frac{r}{p^-_i}}
             |A_{k_{n+1}}|^{1-\frac{r}{(p^-_i)^*}}, \\
       \end{split}
    \end{align}
    where $C=\max(1,C(p^-_1,N), \ldots, C(p^-_m,N))$ with $C(p^-_j,N)$ being the embedding constant
    corresponding to the embedding $W^{1,p^-_j}(\Om) \hookrightarrow
    L^{(p^-_j)^*}(\Om)$, $j=1,\ldots,m$. Thus $C$ is independent of
    $i$. A simple calculation shows that the right-hand side of (\ref{3.13})
can be estimated further to obtain
    \begin{align}
        \begin{split}\label{3.14}
                    & \into (u-k_{n+1})_+^{r}\xi_i^{r}dx\\
                & \leq d_5 \left ( \into |\nabla (u-k_{n+1})_+|^{p^-_i} \xi_i^{p^-_i}dx \right)^{\frac{r}{p^-_i}}|A_{k_{n+1}}|^{1-\frac{r}{(p^-_i)^*}} \\
            & \qquad +d_6 \left (\int_{A_{k_{n+1}}} u^{q_0(x)}dx \right
            )^{\frac{r}{p^-_i}}|A_{k_{n+1}}|^{1-\frac{r}{(p^-_i)^*}}.
        \end{split}
    \end{align}
Here $d_5=d_5(q^+,C)$ and $d_6=d_6(p^+,q^+,C,L)$ are positive
constants, where $L$ is the constant introduced in step (ii).

Using (\ref{3.5}), (\ref{3.11}), (\ref{3.12}), and (\ref{3.14}) with
$r=q^+_{0,i}$ and $r=q^-_{0,i}$, respectively, for $i=1,\ldots,m$,
we get
    \begin{align}\label{3.15}
        \begin{split}
            Z_{n+1}
            & \leq m^{q_0^+} \sum_{i=1}^m \left [
                d_5 \left ( d_4 a^n(Z_n+\tilde{Z}_n)\right)^{\frac{q^+_{0,i}}{p^-_i}}|A_{k_{n+1}}|^{1-\frac{q^+_{0,i}}{(p^-_i)^*}} \right. \\
            & \qquad \qquad
                +d_6 \left (2^{q_0^+(n+2)}Z_n\right )^{\frac{q^+_{0,i}}{p^-_i}}|A_{k_{n+1}}|^{1-\frac{q^+_{0,i}}{(p^-_i)^*}} \\
            & \qquad \qquad
                + d_5 \left ( d_4 a^n(Z_n+\tilde{Z}_n)dx \right)^{\frac{q^-_{0,i}}{p^-_i}}|A_{k_{n+1}}|^{1-\frac{q^-_{0,i}}{(p^-_i)^*}} \\
            & \qquad \qquad  \left.
                +d_6 \left (2^{q_0^+(n+2)}Z_n\right )^{\frac{q^-_{0,i}}{p^-_i}}|A_{k_{n+1}}|^{1-\frac{q^-_{0,i}}{(p^-_i)^*}} \right ].
        \end{split}
    \end{align}
    Setting
    \[
    Y_n:=Z_n+\tilde{Z}_n
    \]
and $\eta=\max\left(\frac{q^+_{0,1}}{(p^-_1)^*}, \ldots,
\frac{q^+_{0,m}}{(p^-_m)^*} \right)$ we have for $r \in
\{q^+_{0,i},q^-_{0,i}\}$
    \begin{align} \label{3.16}
        \begin{split}
            & \left ( d_4 a^n(Z_n+\tilde{Z}_n)\right)^{\frac{r}{p^-_i}} \le
            \left ( d_4 \right )^{\frac{q^+}{p^-}} \left(a^{\frac{q^+}{p^-}}\right )^n (Y_n+Y_n^{\frac{q^+}{p^-}}), \\
            & \left (2^{q^+_0(n+2)}Z_n\right )^{\frac{r}{p^-_i}} \leq \left ( 2^{\frac{(q^+)^2(n+2)}{p^-}}\right) (Y_n+Y_n^{\frac{q^+}{p^-}}), \\
            & |A_{k_{n+1}}|^{1-\frac{r}{(p^-_i)^*}} \leq 2^{q^+(n+1)} \left (\frac{1}{k^{q_0^-}} \right)^{1-\eta}  (Y_n+Y_n^{1-\eta}).
        \end{split}
    \end{align}
    Using these estimates, we conclude from (\ref{3.15}) that
    \begin{align}\label{3.17}
        \begin{split}
            & Z_{n+1}\leq d_7 d_8^n \frac{1}{k^{q_0^- (1-\eta)}} \left
            (Y_n^2+Y_n^{2-\eta}+Y_n^{1+\frac{q^+}{p^-}}+Y_n^{1+\frac{q^+}{p^-}-\eta} \right )
        \end{split}
    \end{align}
    where $d_7$ and $d_8$ are positive constants that only depend on the
    data.

${}$

{\bf (v) Estimating $\tilde{Z}_{n+1}$.} We proceed similarly as in
step (iv). Analogously to (\ref{3.12}) we get an estimate for
$\tilde{Z}_{n+1}$ of the form
    \begin{align}\label{3.18}
        \begin{split}
            \tilde{Z}_{n+1}
            \leq m^{q_1^+} \sum_{i=1}^m \left [\int_{\Gamma_{k_{n+1}}}\!\!\!\!\!\!\! (u-k_{n+1})^{q_{1,i}^+}\xi_i^{q_{1,i}^+}d\sigma +
            \int_{\Gamma_{k_{n+1}}}\!\!\!\!\!\!\! (u-k_{n+1})^{q_{1,i}^-}\xi_i^{q_{1,i}^-} d\sigma \right ],
        \end{split}
    \end{align}
    where $q_{1,i}^-=\min_{B_i\cap \Gamma}q_1(x)$. Note that
$p^-_i \leq q^-_{1,i} \leq q^+_{1,i}<(p^-_i)_*$ for $i=1,\ldots,m$.

Let now $i\in \{1,\ldots,m\}$ be fixed, and suppose that $r\in
\{q_{1,i}^-,q_{1,i}^+\}$, that is we have $p^-_i \leq r <(p^-_i)_*$
and $r \leq q^+$. Define $s=s_i(r)\in (1,N)$ by means of
\begin{align*}
        s_*=
        \begin{cases}
             \frac{r+(p_i^-)_*}{2}\quad & \text{ if } (p_i^-)_* <\infty, \\
            r+1 & \text{ if } (p_i^-)_* =\infty.
        \end{cases}
\end{align*}
Then $s<p^-_i\le r<s_*<(p^-_i)_*$. Since the trace operator maps
$W^{1,s}(\Om)$ boundedly into $L^r(\Gamma)$, we have, using
H\"{o}lder's inequality,
    \begin{align*}
        \begin{split}
            & \int_{\Gamma} ((u-k_{n+1})_+ \xi_i)^{r}d\sigma \\
            & \leq \hat{C}^r \left [\into \left (|\nabla [(u-k_{n+1})_+\xi_i]|^{s}+
            |(u-k_{n+1})_+\xi_i|^{s} \right )dx \right ]^{\frac{r}{s}} \\
            & \leq \hat{C}^{q^+} \left [\into \left (|\nabla [(u-k_{n+1})_+\xi_i]|^{p^-_i}+
            |(u-k_{n+1})_+\xi_i|^{p^-_i} \right )dx \right ]^{\frac{r}{p^-_i}}
            |A_{k_{n+1}}|^{\left (1-\frac{s}{p^-_i} \right ) \frac{r}{s}},
        \end{split}
    \end{align*}
    where $\hat{C}$ denotes the maximum of $1$ and the norms of the trace maps $\gamma: W^{1,s}(\Om)\rightarrow
    L^r(\Gamma)$ when $r$ runs through set
    $\bigcup_{j=1}^m\{q_{1,j}^-,q_{1,j}^+\}$.
    The right-hand side of the last inequality can be estimated to get
    \begin{align}\label{3.19}
        \begin{split}
            & \int_{\Gamma} ((u-k_{n+1})_+ \xi_i)^{r}d\sigma \\
            & \leq d_9 \left ( \into |\nabla (u-k_{n+1})_+|^{p^-_i} \xi_i^{p^-_i}dx \right)^{\frac{r}{p^-_i}}
            |A_{k_{n+1}}|^{\left (1-\frac{s}{p^-_i} \right )\frac{r}{s}} \\
            & \qquad +d_{10} \left (\int_{A_{k_{n+1}}} u^{q_0(x)}dx \right )^{\frac{r}{p^-_i}}
            |A_{k_{n+1}}|^{\left (1-\frac{s}{p^-_i} \right )\frac{r}{s}}
        \end{split}
    \end{align}
    with positive constants $d_9, d_{10}$ that only depend on the
    data.

    From (\ref{3.5}), (\ref{3.11}), (\ref{3.18}), and (\ref{3.19}) with $r=q^+_{1,i}$ and $r=q^-_{1,i}$, respectively, for $i=1,\ldots,m$,
    we infer that
        \begin{align}\label{3.20}
        \begin{split}
            \tilde{Z}_{n+1}
            & \leq m^{q_1^+} \sum_{i=1}^m \left [
                d_9 \left ( d_4a^n(Z_n+\tilde{Z}_n) \right)^{\frac{q^+_{1,i}}{p^-_i}}
                |A_{k_{n+1}}|^{\left (1-\frac{s_i(q^+_{1,i})}{p^-_i} \right )\frac{q^+_{1,i}}{s_i(q^+_{1,i})}} \right. \\
            & \qquad \qquad
                +d_{10} \left (2^{q^+_0(n+2)}Z_n \right )^{\frac{q^+_{1,i}}{p^-_i}}
                |A_{k_{n+1}}|^{\left (1-\frac{s_i(q^+_{1,i})}{p^-_i} \right )\frac{q^+_{1,i}}{s_i(q^+_{1,i})}}
             \\
            & \qquad \qquad
                + d_9 \left ( d_4a^n(Z_n+\tilde{Z}_n) \right)^{\frac{q^-_{1,i}}{p^-_i}}
                |A_{k_{n+1}}|^{\left (1-\frac{s_i(q^-_{1,i})}{p^-_i} \right )\frac{q^-_{1,i}}{s_i(q^-_{1,i})}}
 \\
            & \qquad \qquad  \left.
                + d_{10} \left ((2^{q^+_0(n+2)}Z_n \right )^{\frac{q^-_{1,i}}{p^-_i}}
                |A_{k_{n+1}}|^{\left (1-\frac{s_i(q^-_{1,i})}{p^-_i} \right )\frac{q^-_{1,i}}{s_i(q^-_{1,i})}}
 \right ].
        \end{split}
    \end{align}
Put $\tilde{\eta}=\max\left(\frac{s_1(q^+_{1,1})}{p^-_1}, \ldots,
\frac{s_m(q^+_{1,m})}{p^-_m} \right)$. Similarly to (\ref{3.16})
    we have for $r \in \{q^+_{1,i},q^-_{1,i}\}$
    \begin{align}\label{3.21}
        \begin{split}
            & \left ( d_4 a^n(Z_n+\tilde{Z}_n)\right)^{\frac{r}{p^-_i}} \le \left ( d_4 \right )^{\frac{q^+}{p^-}} \left(a^{\frac{q^+}{p^-}}\right )^n (Y_n+Y_n^{\frac{q^+}{p^-}}), \\
            & \left (2^{q_0^+(n+2)}Z_n\right )^{\frac{r}{p^-_i}} \leq \left ( 2^{\frac{(q^+)^2(n+2)}{p^-}}\right) (Y_n+Y_n^{\frac{q^+}{p^-}}), \\
            & |A_{k_{n+1}}|^{\left (1-\frac{s_i(r)}{p^-_i} \right )\frac{r}{s_i(r)}} \leq 2^{q^+(n+1)}
            \left (\frac{1}{k^{q_0^-}} \right)^{1-\tilde{\eta}}  (Y_n^{q^+}+Y_n^{1-\tilde{\eta}}).
        \end{split}
    \end{align}
    Finally, (\ref{3.20}) and (\ref{3.21}) imply that
    \begin{align}\label{3.22}
        \begin{split}
            & \tilde{Z}_{n+1}\leq d_{11} d_{12}^n \frac{1}{k^{q_0^- (1-\tilde{\eta})}}
            \left (Y_n^{q^+1}+Y_n^{2-\tilde{\eta}}+Y_n^{q^++\frac{q^+}{p^-}}+Y_n^{1+\frac{q^+}{p^-}-\tilde{\eta}} \right )
        \end{split}
    \end{align}
     where $d_{11}$ and $d_{12}$ are positive constants that only depend on the
    data.

${}$

{\bf (vi) The iterative inequality for $Y_{n}$.} Recall that
$Y_n=Z_n+\tilde{Z}_n$. Hence (\ref{3.17}) and (\ref{3.22}) yield
    \begin{align*}
        \begin{split}
        Y_{n+1}
        & \leq K b^n \frac{1}{k^{q_0^- (1-\hat{\eta})}}
         \left ( Y_n^2+Y_n^{2-\eta}+Y_n^{1+\frac{q^+}{p^-}}+Y_n^{1+\frac{q^+}{p^-}-\eta} \right.\\
        & \quad\quad \left. \quad +Y_n^{q^+}+Y_n^{2-\tilde{\eta}}+Y_n^{q^++\frac{q^+}{p^-}}+Y_n^{1+\frac{q^+}{p^-}-\tilde{\eta}} \right
        )\\
        & \leq 8K b^n \frac{1}{k^{q_0^-
        (1-\hat{\eta})}}\big(Y_n^{1+\delta_1}+Y_n^{1+\delta_2}\big)
        \end{split}
    \end{align*}
    with $K=\max(d_{7},d_{11}), b=\max(d_8,d_{12})$,
    $\hat{\eta}=\max(\eta,\tilde{\eta})$, and where $0<\delta_1\le
    \delta_2$ are given by
\begin{align*}
\delta_1=&\;\min(1,1-\eta,\frac{q^+}{p^-},\frac{q^+}{p^-}-\eta,q^+,1-\tilde{\eta},q^++\frac{q^+}{p^-}-1,\frac{q^+}{p^-}-\tilde{\eta}),\\
\delta_2=&\;\max(1,1-\eta,\frac{q^+}{p^-},\frac{q^+}{p^-}-\eta,q^+,1-\tilde{\eta},q^++\frac{q^+}{p^-}-1,\frac{q^+}{p^-}-\tilde{\eta}).
    \end{align*}
    Without loss of generality we may assume that $b>1$.
Now we may apply Lemma \ref{lemma4}, which says that $Y_n \ra 0$ as
$n\to \infty$ provided
    \begin{align}\label{3.23}
        Y_0=\into (u-k)_+^{q_0(x)} dx +\int_{\Gamma} (u-k)_+^{q_1(x)} d \sigma \leq
        \left (\frac{16 K}{k^{q_0^- (1-\hat{\eta})}} \right )^{-\frac{1}{\delta_1}}
        b^{-\frac{1}{\delta_1^2}}.
    \end{align}
Relation (\ref{3.23}) is clearly satisfied if
    \begin{align}\label{3.24}
        \into u_+^{q_0(x)} dx + \int_{\Gamma} u_+^{q_1(x)} d \sigma \int_{\Gamma}\leq
        \left (\frac{16 K}{k^{q_0^- (1-\hat{\eta})}} \right )^{-\frac{1}{\delta_1}} b^{-\frac{1}{\delta_1^2}}.
    \end{align}
    Hence, if we choose $k$ such that
    \begin{align}
        k =\max\left(1, \left [ (16K)^{\frac{1}{\delta_1}} b^{\frac{1}{\delta_1^2}}
        \left (\into u_+^{q_0(x)} dx + \int_{\Gamma} u_+^{q_1(x)} d \sigma \right )\right ]^{\frac{\delta_1}{q^-_0(1-\hat{\eta})}}\right),
    \end{align}
    then (\ref{3.24}) and in particular (\ref{3.23}) are satisfied. Since $k_n \ra 2k$ as $n \ra \infty$ we obtain
    \begin{align*}
        \esssup_{\Om} u \leq 2k = 2
        \max\left(1, \left [ (16K)^{\frac{1}{\delta_1}} b^{\frac{1}{\delta_1^2}} \left (\into u_+^{q_0(x)} dx +
         \int_{\Gamma} u_+^{q_1(x)} d \sigma \right )\right
         ]^{\frac{\delta_1}{q^-_0(1-\hat{\eta})}}\right).
    \end{align*}
Tracing back the constants, we see that the first part of the
theorem is proved. The supersolution case can be done analogously,
replacing $u$ with $-u$ and $A_k$ with $\tilde{A}_k$, and using
Lemma \ref{lemma2} instead of Lemma \ref{lemma1}. This completes the
proof.
\end{proof}
%

\end{document}


\maketitle
In this corrigendum we correct a lemma concerning the geometric convergence of sequences of numbers which was used in \cite{Winkert-Zacher-2012} as Lemma 2.1. As a consequence the statement in the main result changes a bit and the corresponding proof needs some minor different arguments to be fitted.\\
{\bf (a)} First, we replace Theorem 1.1 in \cite{Winkert-Zacher-2012} by the following one:
\setcounter{section}{1}
\setcounter{theorem}{0}
\begin{theorem}\label{maintheorem}
    Let the assumptions in (H) be satisfied. Then there exist
    positive constants $\alpha=\alpha(p,q_0,q_1)$ and
    $C=C(p,q_0,q_1,a_3,a_4,a_5,b_0,b_1,b_2,c_0,c_1,N,\Omega)$ such
    that the following assertions hold.
    \begin{enumerate}
        \item[(i)]
            If $u \in \Wx1p$ is a weak subsolution of (1.1), then both $\esssup_{\Om} u $ and $\esssup_{\Gamma}u$ are bounded from above by 
            \begin{align*}
                C \left [1+\into u_+^{q_0(x)} dx + \int_{\Gamma} u_+^{q_1(x)} d \sigma \right]^\alpha.
            \end{align*}
        \item[(ii)]
            If $u \in \Wx1p$ is a weak supersolution of (1.1), then both $\essinf_{\Om}u$ and $\essinf_{\Gamma}u$ are bounded from below by 
            \begin{align*}
                -C\left [1+\into (-u)_+^{q_0(x)}  dx + \int_{\Gamma} (-u)_+^{q_1(x)}   d \sigma \right ]^\alpha.
            \end{align*}
    \end{enumerate}
\end{theorem}
\noindent{\bf (b)}
Next, we replace Corollary 1.2 in \cite{Winkert-Zacher-2012} by the following one: 
\begin{corollary}
    Let the assumptions (H) be satisfied and let $u \in \Wx1p$ be a weak solution of (1.1).  Then $u\in L^\infty(\Omega), L^\infty(\Gamma)$ and the estimates in (i) and (ii)
    from Theorem \ref{maintheorem} are valid.
\end{corollary}
\noindent{\bf (c)} Replace reference [32] on page 4, line 5 from bottom by the new reference \cite{Ho-Sim-2015}.\\
{\bf (d)} Now, we replace Lemma 2.1 in \cite{Winkert-Zacher-2012} by the following one:
\setcounter{section}{2}
\setcounter{theorem}{0}
\begin{lemma}\label{lemma_geometric_convergence}
    Let $\{Y_n\}, n=0,1,2,\ldots,$ be a sequence of positive numbers, satisfying the recursion inequality
    \begin{align*}
        Y_{n+1} \leq K b^n \left (Y_n^{1+\delta_1}+ Y_n^{1+\delta_2} \right ) , \quad n=0,1,2, \ldots,
    \end{align*}
    for some $b>1,\,K>0$ and $\delta_2\geq \delta_1>0$. If
    \begin{align*}
        Y_0 \leq \min \left(1,(2K)^{-\frac{1}{\delta_1}} b^{-\frac{1}{\delta_1^2}}\right)
    \end{align*}
    or
    \begin{align*}
        Y_0 \leq \min  \left((2K)^{-\frac{1}{\delta_1}} b^{-\frac{1}{\delta_1^2}}, (2K)^{-\frac{1}{\delta_2}}b^{-\frac{1}{\delta_1 \delta_2}-\frac{\delta_2-\delta_1}{\delta_2^2}}\right),
    \end{align*}
    then $Y_n \leq 1$ for some $n \in \N \cup \{0\}$. Moreover,
    \begin{align*}
        Y_n \leq \min \left(1,(2K)^{-\frac{1}{\delta_1}} b^{-\frac{1}{\delta_1^2}} b^{-\frac{n}{\delta_1}}\right), \quad \text{ for all }n \geq n_0,
    \end{align*}
    where $n_0$ is the smallest $n \in \N \cup \{0\}$ satisfying $Y_n \leq 1$. In particular, $Y_n \to 0$ as $n \to \infty$.
\end{lemma}

We note that Lemma 2.1 stated in \cite{Winkert-Zacher-2012} would have been correct if $K>1$ instead of $K>0$. However, we need in our treatment such a result for arbitrary positive $K$. 

Now, at two places in the proof of Theorem 1.1, we need some minor changes.

\noindent{\bf (e)} On page 8, after line 3, we add the following paragraph:
	
``Here, $\left(p_i^-\right)^*$ and $\left(p^-_i\right)_*$ are defined by, for all $i=1,\ldots,m$, 
\begin{align*}
    \left(p_i^-\right)^*=
    \begin{cases}
	\frac{N\left(p_i^-\right)}{N-\left(p_i^-\right)} & \text{ if } \left(p_i^-\right) <N, \\
	q_{0}^++1 & \text{ if } \left(p_i^-\right) \geq N,
    \end{cases} \quad
    \left(p^-_i\right)_*=
    \begin{cases}
	\frac{(N-1)\left(p^-_i\right)}{N-\left(p^-_i\right)} & \text{ if } \left(p^-_i\right) <N, \\
	q_{1}^++1 & \text{ if } \left(p^-_i\right) \geq N,
    \end{cases}
\end{align*}
where $q^+_0=\max_{x\in\overline{\Omega}}q_0(x)$ and $q^+_1=\max_{x\in\Gamma}q_1(x)$ (see Section 2).''\\

\noindent{\bf (f)} Replace the paragraph on page 12 from formula (3.23) until line 4 from bottom by the following paragraph:\\
\setcounter{section}{3}\setcounter{equation}{22}
``
\begin{align}\label{VI_1}
    \begin{split}
	Y_0
	& =\into (u-k)_+^{q_0(x)} dx +\int_{\Gamma} (u-k)_+^{q_1(x)} d \sigma \\
	& \leq \min \left[\left (\frac{16 K}{k^{q_0^- (1-\hat{\eta})}} \right )^{-\frac{1}{\delta_1}} b^{-\frac{1}{\delta_1^2}}, \left (\frac{16 K}{k^{q_0^- (1-\hat{\eta})}} \right )^{-\frac{1}{\delta_2}}b^{-\frac{1}{\delta_1 \delta_2}-\frac{\delta_2-\delta_1}{\delta_2^2}}\right].
    \end{split}
\end{align}
Relation (\ref{VI_1}) is clearly satisfied if
\begin{align}\label{VI_2}
    \begin{split}
	&\into u_+^{q_0(x)} dx + \int_{\Gamma} u_+^{q_1(x)} d \sigma\\
	& \leq \min \left[\left (\frac{16 K}{k^{q_0^- (1-\hat{\eta})}} \right )^{-\frac{1}{\delta_1}} b^{-\frac{1}{\delta_1^2}}, \left (\frac{16 K}{k^{q_0^- (1-\hat{\eta})}} \right )^{-\frac{1}{\delta_2}}b^{-\frac{1}{\delta_1 \delta_2}-\frac{\delta_2-\delta_1}{\delta_2^2}}\right].
    \end{split}
\end{align}
Hence, if we choose $k$ such that
\begin{align}\label{VI_3}
    \begin{split}
	k
	& =\left(1+\vphantom{\left[ \int_0^\delta\int_{\Gamma} u_+^{q_2(t,x)} d \sigma dt\right]^{\frac{\delta_1}{q^-_1(1-\hat{\eta})}}} 
	(16K)^{\frac{1}{q_0^- (1-\hat{\eta})}} b^{\frac{1}{\delta_1 q_0^- (1-\hat{\eta})}+\frac{\delta_2-\delta_1}{\delta_2 q_0^- (1-\hat{\eta})}}\right)\\
	& \qquad\qquad  \times 
	\left (1+\into u_+^{q_0(x)} dx + \int_{\Gamma} u_+^{q_1(x)} d \sigma \right )^{\frac{\delta_2}{q^-_0(1-\hat{\eta})}},
    \end{split}
\end{align}
then (\ref{VI_2}) and in particular (\ref{VI_1}) are satisfied. Since $k_n \ra 2k$ as $n \ra \infty$ we obtain
\begin{align*}
    \esssup_{\Om} u \leq 2k \quad \text{and} \quad \esssup_{\Gamma} u \leq 2k
\end{align*}
with $k$ given in \eqref{VI_3}.
''